\documentclass[11pt]{article}
\usepackage{amssymb,amsfonts}
\usepackage{amsthm}
\usepackage{graphicx}
\newcommand{\overbar}[1]{\mkern 3mu\overline{\mkern-3mu#1\mkern-0mu}\mkern 0mu}

\newcommand{\twomat}[4]{\left(\begin{array}{cc}#1&#2\\#3&#4\end{array}\right)}

\newcommand{\R}{{\mathbb{R}}}

\DeclareRobustCommand\openone{\leavevmode\hbox{\small1\normalsize\kern-.33em1}}

\newcommand{\be}{\begin{equation}}
\newcommand{\ee}{\end{equation}}
\newcommand{\bea}{\begin{eqnarray}}
\newcommand{\eea}{\end{eqnarray}}
\newcommand{\beas}{\begin{eqnarray*}}
\newcommand{\eeas}{\end{eqnarray*}}

\newtheorem{theorem}{Theorem}

\newcount\minute
\newcount\hour
\def\currenttime{%
    \minute\time
    \hour\minute
    \divide\hour60
    \the\hour:\multiply\hour60\advance\minute-\hour\the\minute}
\begin{document}
\title{A generalisation of Mirsky's singular value inequalities}
\author{Koenraad M.R.\ Audenaert\\[3mm]
Department of Mathematics,
Royal Holloway University of London, \\
Egham TW20 0EX, United Kingdom \\[1mm]
Department of Physics and Astronomy, University of Ghent, \\
S9, Krijgslaan 281, B-9000 Ghent, Belgium\\[2mm]
koenraad.audenaert@rhul.ac.uk}
\date{\today, \currenttime}
\maketitle
\begin{abstract}
We prove an $f$-version of Mirsky's singular value inequalities for differences of matrices. This $f$-version
consists in applying a positive concave function $f$, with $f(0)=0$, 
to every singular value in the original Mirsky inequalities.
\end{abstract}

Denote the singular values of a matrix $X$, arranged in non-increasing order, by $\sigma_i(X)$.
The main result of this paper is the following singular value inequality:
\begin{theorem}\label{th:main}
Let $f$ be a concave function $f:\R_+\mapsto\R_+$ with $f(0)=0$.
Let $X,Y$ be general $n\times n$ complex matrices.
Then for any $m\le n$, and any increasing sequence $(i_1,\ldots,i_m)$ of integers in $\{1,\ldots,n\}$, we have
\be
\sum_{k=1}^m |f(\sigma_{i_k}(X)) - f(\sigma_{i_k}(Y))| \le
\sum_{k=1}^m f(\sigma_k(X-Y)).\label{eq:0}
\ee
\end{theorem}
Without the application of the function $f$, these inequalities are essentially Mirsky's singular value inequalities \cite{M} 
(up to setting $i_k=k$).
Their $f$-version is essentially the set of inequalities conjectured by W.~Miao that appears in \cite{AK} as Conjecture 6 
(again with $i_k=k$).
We therefore have solved this conjecture.
The special case $i_k=k$ and $m=n$ has apparently been proven by Yue and So in their as yet unpublished manuscript \cite{YS},
where an application is given to low-rank matrix recovery.
Our proof technique is completely different from theirs.
In \cite{ZQ} Zhang and Qiu also claimed to have proven inequalities (\ref{eq:0}), 
but unfortunately their proof is flawed (as pointed out in \cite{YS}).

The main ingredient in our work is a set of eigenvalue inequalities for sums of Hermitian matrices,
known as the Thompson-Freede (TF) inequalities \cite{TF}. These inequalities also come in a version that applies
to singular values.
Remarkably, this is about the \emph{only} matrix analytical tool
that is needed to prove Theorem \ref{th:main}. Apart from this, the proof is rather elementary and consists in 
appropriately choosing one of the TF inequalities and combining it with inequalities of the kind 
$\sigma_{i}(X)\ge\sigma_j(X)$ for $i<j$, and $\sigma_i(X)\ge0$.

In Section \ref{sec:1} we introduce the TF inequalities, for eigenvalues as well as for singular values.
We then state
their $f$-version, Theorem \ref{th:eig}, by which is meant applying a positive, concave function $f$ with $f(0)=0$ 
to every singular value in
the original TF singular value inequalities. Zhang and Qiu have shown in \cite{ZQ} that all Horn-type singular value 
inequalities have a valid $f$ version, including the TF inequalities.
We give a completely different proof of the $f$-version of the TF singular value inequalities 
that just like the proof of Theorem \ref{th:main}
is only based on a well-chosen combination of the original TF inequalities. In fact, neither the statement of the theorem
nor its proof make any reference to matrix analysis at all.
The proof of our main result, Theorem \ref{th:main}, is given in Section \ref{sec:2}.
\section{The Thompson-Freede inequalities and their $f$-version\label{sec:1}}

Let $A$ and $B$ be $n\times n$ Hermitian matrices.
Let $\alpha(i)$, $\beta(i)$ and $\gamma(i)$, for $i=1,\ldots,n$, be the eigenvalues,
sorted in non-ascending order, of $A$, $B$ and $A+B$, respectively.

For a given integer $m\le n$ let $(i_1,\ldots,i_{m})$ and $(j_1,\ldots,j_{m})$ be
two strictly increasing sequences of length $m$ of integers between 1 and $n$
such that $i_{m}+j_{m}\le n+m$.
Then the Thompson-Freede (TF) inequalities \cite{TF} are
\be
\sum_{k=1}^{m} \gamma(i_k+j_k-k) \le \sum_{k=1}^{m} \alpha(i_k)+\sum_{k=1}^{m} \beta(j_k).\label{eq:TF}
\ee
These inequalities include as special cases the Lidskii/Wielandt inequalities (take $j_k=k$).
They form themselves a subset of Horn's inequalities, which are of the general form
\be
\sum_{k\in K} \gamma(k) \le \sum_{i\in I} \alpha(i)+\sum_{j\in J} \beta(j),\label{eq:horn}
\ee
where $I$, $J$ and $K$ are certain subsets of $\{1,\ldots,n\}$ governed by a rather 
complicated set of recursive constraints (not reproduced here) \cite{bhat}.
We say that $(I,J,K)$ constitutes an \emph{admissible triple} whenever these constraints are satisfied.

Consider now a non-negative, concave (hence non-decreasing) function $f$ on $[0,+\infty)$ such that $f(0)=0$.
When $A$ and $B$ are positive semidefinite, the eigenvalues of $A$, $B$ and $A+B$ also 
satisfy what one could call the $f$-version of Horn's inequalities:
\be
\sum_{k\in K} f(\gamma(k)) \le \sum_{i\in I} f(\alpha(i))+\sum_{j\in J} f(\beta(j)).\label{eq:fhorn}
\ee
These inequalities are also satisfied for \emph{general} matrices $A$ and $B$
when $\alpha$, $\beta$ and $\gamma$ are the \emph{singular values} of $A$,
$B$ and $A+B$, respectively.
Note that the non-negativity of singular values is essential here;
although the eigenvalues of Hermitian $A$, $B$ and $A+B$ satisfy all Horn inequalities,
they do not in general satisfy their $f$-versions.

Zhang and Qiu \cite{ZQ} have recently proven this $f$-version by exploiting a theorem by Bourin and Uchiyama 
(Corollary 2.6 in \cite{BU}) which states that
for all $A$ and $B$ and any positive concave function $f$ there exist unitary matrices 
$U$ and $V$ such that
\be
f(|A+B|) \le U f(|A|) U^* + V f(|B|) V^*.\label{eq:BU}
\ee
Thus, in particular, the singular values of $A$, $B$ and $A+B$ satisfy 
an $f$-version of the TF inequalities.

Below we give an alternative proof of the latter statement based \emph{uniquely} on the fact that 
these singular values satisfy the original TF inequalities.
\begin{theorem}\label{th:eig}
Let $\alpha(i)$, $\beta(i)$ and $\gamma(i)$, for $i=1,\ldots,n$, be sequences of non-negative numbers,
sorted in non-ascending order, and satisfying all TF inequalities (\ref{eq:TF}).
Let $f$ be a non-negative, concave function on $[0,+\infty)$ such that $f(0)=0$.
Then $\alpha(i)$, $\beta(i)$ and $\gamma(i)$ satisfy the $f$-version of the TF inequalities.
To wit, for a given integer $m\le n$ let $(i_1,\ldots,i_{m})$ and $(j_1,\ldots,j_{m})$ be
two strictly increasing sequences of length $m$ of integers between 1 and $n$
such that $i_{m}+j_{m}\le n+m$.
Then
\be
\sum_{k=1}^m f(\gamma(i_k+j_k-k)) \le \sum_{k=1}^m f(\alpha(i_k))+\sum_{k=1}^m f(\beta(j_k)).\label{eq:eig}
\ee
\end{theorem}
\textit{Proof.}
Any function $f$ satisfying the assumptions of Theorem \ref{th:eig} can be uniformly
approximated as a finite or infinite positive linear combination of `hook' functions
$h_t(x):=\min(x,t)$ with $t>0$; that is, for any such $f$ there exists a positive
measure $d\mu(t)$ on $(0,\infty)$ such that
$f(x)=\int_0^\infty h_t(x)\;d\mu(t)$. By linearity of LHS and RHS of (\ref{eq:eig}) in $f$
it therefore suffices to prove (\ref{eq:eig}) for hook functions only. Furthermore,
by a scaling argument it is clear that we can restrict to $f(x)=h(x):=\min(x,1)$.

Let $a$ and $b$ be index values, $1\le a,b\le m$, such that the following hold:
\[
\alpha(i_a) < 1\le \alpha(i_{a-1}), \quad
\beta(j_b)  < 1\le \beta(j_{b-1}).
\]
Then
$h(\alpha(i_k))=1$ for $k<a$ and $h(\alpha(i_k))=\alpha(i_k)$ for $k\ge a$,
and similar identities hold for $\beta$.
Inequality (\ref{eq:eig}) then reduces to
\be
\sum_{k=1}^m h(\gamma(i_k+j_k-k)) \le (a-1)+\sum_{k=a}^m \alpha(i_k) + (b-1) + \sum_{k=b}^m \beta(j_k).\label{eq:eig2}
\ee

We will first consider the case that $a+b-1\le m$.
Making the replacements $m\to m':=m-(a-1)-(b-1)$, 
$i_k \to i_{k+a-1}$ and $j_k \to j_{k+b-1}$ in (\ref{eq:TF}) yields
\bea
\sum_{k=1}^{m'} \gamma(i_{k+a-1} + j_{k+b-1}-k)
&\le& \sum_{k=1}^{m'} \alpha(i_{k+a-1}) + \sum_{k=1}^{m'} \beta(j_{k+b-1}) \nonumber\\
&=& \sum_{k=a}^{m-b+1} \alpha(i_{k}) + \sum_{k=b}^{m-a+1} \beta(j_{k}).\label{eq:eig2c}
\eea
Since $(i_1,\ldots,i_m)$ and $(j_1,\ldots,j_m)$ are strictly increasing sequences, they satisfy
$i_{k+a-1}\le i_{k+a+b-2}-(b-1)$ and $j_{k+b-1}\le j_{k+a+b-2}-(a-1)$.
Furthermore, the sequence $(\gamma(1),\ldots,\gamma(n))$ is non-increasing. Therefore,
the LHS of (\ref{eq:eig2c}) is bounded below as
\beas
\sum_{k=1}^{m'} \gamma(i_{k+a-1}+j_{k+b-1}-k) &\ge& \sum_{k=1}^{m'} \gamma(i_{k+a+b-2}+j_{k+a+b-2}-(k+a+b-2))\\
 &=& \sum_{k=a+b-1}^{m} \gamma(i_{k}+j_{k}-k),
\eeas
so that
\[
\sum_{k=a+b-1}^m \gamma(i_k+j_k-k) \le \sum_{k=a}^{m-b+1} \alpha(i_k) + \sum_{k=b}^{m-a+1} \beta(j_k).
\]

For the other case, $a+b-1>m$, the same inequality holds trivially.

Because $h(x)\le x$ and $\alpha,\beta\ge0$ (here is where the argument would break down 
when considering eigenvalues instead of singular values) this implies
\be
\sum_{k=a+b-1}^m h(\gamma(i_k+j_k-k)) \le \sum_{k=a+b-1}^m \gamma(i_k+j_k-k)
\le \sum_{k=a}^m \alpha(i_k) + \sum_{k=b}^m \beta(k).\label{eq:zz1}
\ee
For the remaining terms in the LHS of (\ref{eq:eig2}) we have
\be
\sum_{k=1}^{a+b-2} h(\gamma(i_k+j_k-k))\le a+b-2 = (a-1)+(b-1),\label{eq:zz2}
\ee
and taking the sum of (\ref{eq:zz1}) and (\ref{eq:zz2}), inequality (\ref{eq:eig2}) follows.
\qed
\section{Proof of Theorem \ref{th:main}\label{sec:2}}
Let us replace the matrices $X$, $Y$ and $X-Y$ in the statement of Theorem \ref{th:main} 
by matrices $C$, $A$ and $B$, respectively, with $A+B+C=0$,
and let us denote their singular values by $\gamma(i)$, $\alpha(i)$ and $\beta(i)$, respectively.
The proof starts with a number of simple reductions.

As in the proof of Theorem \ref{th:eig}, 
it is enough to prove Theorem \ref{th:main} for the function $f(x):=h(x)=\min(1,x)$, as all other functions under consideration
can be written as positive linear combinations of $t\;h(x/t)$.
Whereas in the proof of Theorem \ref{th:eig} we merely exploited linearity of the LHS and RHS in $f$, 
here we must also use the triangle inequality for the absolute value in the LHS.

It therefore suffices to prove the following inequality
\be
\sum_{k=1}^m |\min(1,\gamma(i_k)) - \min(1,\alpha(i_k))| \le \sum_{k=1}^m \min(1,\beta(k)).\label{eqh:0}
\ee
Let us define the index set $I=\{i_1,\ldots,i_m\}$, and the indices $a$, $b$ and $c$ for which the following holds:
\beas
\gamma(i_c) &< 1\le& \gamma(i_{c-1}) \\
\alpha(i_a) &< 1\le& \alpha(i_{a-1}) \\
\beta(b) &< 1\le& \beta(b-1).
\eeas
We can assume that $a\le c$; otherwise we just swap the roles of $A$ and $C$.

As the contribution to the LHS of (\ref{eqh:0}) of the terms with $k<a$ is exactly zero, 
removing the indices $i_1,\ldots,i_{a-1}$ from $I$ and removing the $a-1$ smallest $\beta$'s from the RHS
turns one instance of (\ref{eqh:0}) into another.
Thus, henceforth we only need to consider the case $a=1$, which is:
\be
\sum_{k=1}^{c-1} (1-\alpha(i_k))
+ \sum_{k=c}^m |\gamma(i_k)-\alpha(i_k)|
\le (b-1) + \sum_{k=b}^m \beta(k).\label{eqh:1}
\ee

Let us partition $I=\{i_1,\ldots,i_m\}$ into two subsets $I_C$ and $I_A$, where $I_C$ is the set of indices 
$i\in I$ for which $\gamma(i)\ge\alpha(i)$ and $I_A$ is the set of
remaining indices. Clearly, the indices $i_1,\ldots,i_{c-1}$ are always in $I_C$, and never in $I_A$.
Because $|x-y| = \max(x-y,y-x)$,
to prove (\ref{eqh:1}) it suffices to prove the following inequality for \emph{all} such partitions $I_C$ and $I_A$ of $I$ 
(keeping the requirement that $i_1,\ldots,i_{c-1}\in I_C$), 
regardless for which of the $i$ the inequality $\gamma(i)\ge\alpha(i)$ holds:
\[
\sum_{k\in I_C} (h(\gamma(k))-\alpha(k))
-\sum_{k\in I_A} (\gamma(k)-\alpha(k))
 \le (b-1) + \sum_{k=b}^m \beta(k),
\]
or, equivalently
\be
\sum_{k\in I_C} h(\gamma(k))+\sum_{k\in I_A} \alpha(k)
 \le \sum_{k\in I_C} \alpha(k) + \sum_{k\in I_A} \gamma(k)
 + (b-1) + \sum_{k=b}^m \beta(k).\label{eqh:3}
\ee

After these reductions, we come to the core of the argument.
Let us define the additional index sets
\bea
I_R&=&\{i_1,\ldots,i_{m-b+1}\},\quad I_{\overbar{R}} = \{i_{m-b+2},\ldots,i_m\}, \nonumber\\
I_L&=&\{i_b,\ldots,i_m\}, \qquad I_{\overbar{L}} = \{i_1,\ldots,i_{b-1}\}, \nonumber\\
J&=&\{b,\ldots,m\}.
\eea
Note that $I_R$, $I_L$ and $J$ have size $m-b+1$ and $I_{\overbar R}$ and $I_{\overbar L}$ have size $b-1$. 
To simplify notations we adopt the notations $\gamma(K):=\sum_{k\in K}\gamma(k)$, etc.,
and $I_{CL}:=I_C\cap I_L$, $I_{C\overbar L}:=I_C\cap I_{\overbar L}$, etc.

Inequality (\ref{eqh:3}), and hence the inequality of Theorem \ref{th:main}, 
is a straightforward consequence of the following theorem, which will be proven below:
\begin{theorem}\label{th:interm}
For all $n\times n$ matrices $A$, $B$ and $C$ such that $A+B+C=0$,
for any partitioning of $I=\{i_1,\ldots,i_m\}$ (with $m\le n$) into $I=I_C\cup I_A$,
and with the notations just introduced, 
\be
\gamma(I_{CL}) + \alpha(I_{AL}) \le \alpha(I_{CR})+\gamma(I_{AR}) + \beta(J).\label{eqh:4}
\ee
\end{theorem}
Note we do not restrict $I_C$ to contain $i_1,\ldots,i_{c-1}$ here.

The simplest non-trivial examples of inequality (\ref{eqh:4}) are 
\be
\gamma(I_{L})  \le \alpha(I_{R}) + \beta(J),\label{eq:ex}
\ee
which are obtained by setting $I_C=I$ and $I_A=\emptyset$. One can easily verify that these are just instances of the TF 
singular value inequalities.
What Theorem \ref{th:interm} is actually saying is that in (\ref{eq:ex}) one can freely replace any $\alpha(i)$ with the
corresponding $\gamma(i)$ and vice-versa, and still have a valid inequality.

\medskip

To see how (\ref{eqh:3}) follows from this, note that
all terms in the LHS of (\ref{eqh:3}) are bounded above by 1.
Therefore, and because $I_{C\overbar L}\cup I_{A\overbar L} = I_{\overbar L}$,
\[
\gamma(I_{C\overbar L}) + \alpha(I_{A\overbar L}) \le b-1.
\]
Furthermore, as singular values are non-negative, we also have
\[
0 \le \alpha(I_{C\overbar R}) + \gamma(I_{A\overbar R}).
\]
Adding these two inequalities to inequality (\ref{eqh:4}) of Theorem \ref{th:interm},
we get
\[
\gamma(I_{C}) + \alpha(I_{A}) \le \alpha(I_{C})+\gamma(I_{A}) + \beta(J).
\]
Since $h(x)=\min(1,x)\le x$ this yields (\ref{eqh:3}).
\qed

\bigskip

\noindent\textit{Proof of Theorem \ref{th:interm}.} 

The essential idea is to consider the following pairing of elements $t_k$ of $I_R$ with elements $s_k$ of $I_L$:
\[
(t_k,s_k):=(i_k,i_{k+b-1}), \quad k=1,\ldots,m-b+1.
\]
Clearly, we have $s_k-t_k =i_{k+b-1}-i_k\ge b-1$.
For every such pair, exactly one out of four possibilities arises
concerning membership of the sets $I_{AL},I_{CL},I_{AR},I_{CR}$.
We can partition the set $K:=\{1,\ldots,m-b+1\}$ accordingly as $K=K_1\cup K_2\cup K_3\cup K_4$, with
\bea
K_1 &=& \{k:s_k\in I_{CL}, t_k\in I_{AR}\} \nonumber\\
K_2 &=& \{k:s_k\in I_{AL}, t_k\in I_{CR}\} \nonumber\\
K_3 &=& \{k:s_k\in I_{CL}, t_k\in I_{CR}\} \nonumber\\
K_4 &=& \{k:s_k\in I_{AL}, t_k\in I_{AR}\}.\label{eq:Ks}
\eea
These sets have the following unions:
\bea
K_1\cup K_3 &=& \{k: s_k \in I_{CL}\} \nonumber\\
K_2\cup K_4 &=& \{k: s_k \in I_{AL}\} \nonumber\\
K_1\cup K_4 &=& \{k: t_k \in I_{AR}\} \nonumber\\
K_2\cup K_3 &=& \{k: t_k \in I_{CR}\}.\label{eq:Kcups}
\eea

With the four subsets $K_i$ as a starting point, we will write down a number of valid inequalities,
the sum of which is exactly (\ref{eqh:4}).

For every $k\in K_1$ we consider the inequality
$\gamma(s_k) \le \gamma(t_k)$, which is valid since $s_k\ge t_k$.
Summing over all $k\in K_1$, we get
\be
\sum_{k\in K_1} \gamma(s_k)
\le
\sum_{k\in K_1} \gamma(t_k).\label{eq:sum1}
\ee
Analogously, we have
\be
\sum_{k\in K_2} \alpha(s_k)
\le
\sum_{k\in K_2} \alpha(t_k).\label{eq:sum2}
\ee

For the remaining pairs, corresponding to $k\in K_3\cup K_4$,
we will write down a single, but more complicated inequality.  
Letting $r$ denote the number of these remaining pairs,
$r=|K_3|+|K_4|$, we write
\be
\sum_{k\in K_3} \gamma(s_k) +\sum_{k\in K_4} \alpha(s_k)
\le
\sum_{k\in K_3} \alpha(t_k) +\sum_{k\in K_4} \gamma(t_k)
+\sum_{k=b}^{b+r-1} \beta(k).\label{eq:sum3}
\ee
Since the overall number of pairs is exactly $m-b+1$, we have $r\le m-b+1$, so that
$b+r-1\le m$. The index $k$ in the final summation therefore does not exceed the bound $m$.

Assuming (\ref{eq:sum3}) is correct, the sum of (\ref{eq:sum1}), (\ref{eq:sum2}) and (\ref{eq:sum3}) yields, after adding 
some more $\beta$-terms to the RHS,
\[
\sum_{k\in K_1\cup K_3} \gamma(s_k) + \sum_{k\in K_2\cup K_4} \alpha(s_k)
\le \sum_{k\in K_1\cup K_4} \gamma(t_k) + \sum_{k\in K_2\cup K_3} \alpha(t_k)
+ \sum_{k=b}^m \beta(k).
\]
By the identities (\ref{eq:Kcups}), this is exactly inequality (\ref{eqh:4}).

It remains to prove inequality (\ref{eq:sum3}). We will do so by showing that 
it is essentially one of the TF inequalities.
Rearranging (\ref{eq:sum3}) gives
\be
\sum_{k\in K_3} \gamma(s_k) -\sum_{k\in K_4} \gamma(t_k)
\le
\sum_{k\in K_3} \alpha(t_k) -\sum_{k\in K_4} \alpha(s_k)
+\sum_{k=b}^{b+r-1} \beta(k).\label{eqh:5}
\ee
The TF inequality that we need is the eigenvalue inequality
\be
\sum_{l=1}^r \widehat\gamma(i_l+j_l-l)
\le \sum_{l=1}^r \widehat\alpha(i_l) + \widehat\beta(j_l)
\label{eq:TFW}
\ee
for eigenvalues $\widehat\alpha(k)$, $\widehat\beta(k)$ and $\widehat\gamma(k)$
of Hermitian matrices $\widehat A$, $\widehat B$ and $\widehat A+\widehat B$, respectively.
In particular, we take $j_l = l+b-1$ (so that $j_l-l=b-1$), and let the indices $i_l$
be the elements of the set
\[
\{t_k: k\in K_3\} \cup \{2n+1-s_k: k\in k_4\}
\]
sorted in decreasing order.
Then (\ref{eq:TFW}) becomes
\beas
\lefteqn{\sum_{k\in K_3} \widehat\gamma(t_k+b-1) + \sum_{k\in K_4} \widehat\gamma(2n+1-s_k+b-1)} \\
&\le& \sum_{k\in K_3} \widehat\alpha(t_k) + \sum_{k\in K_4}\widehat\alpha(2n+1-s_k)
+\sum_{l=1}^r  \widehat\beta(l+b-1).
\eeas
Because $s_k\ge t_k+b-1$, this implies the weaker inequality
\bea
\lefteqn{\sum_{k\in K_3} \widehat\gamma(s_k) + \sum_{k\in K_4} \widehat\gamma(2n+1-t_k)} \nonumber \\
&\le& \sum_{k\in K_3} \widehat\alpha(t_k) + \sum_{k\in K_4}\widehat\alpha(2n+1-s_k)
+\sum_{l=b}^{b+r-1}  \widehat\beta(l).\label{eq:TFW2}
\eea

If we let $\widehat A$ and $\widehat B$ be the Wielandt matrices
\[
\widehat A=\twomat{0}{A}{A^*}{0},\quad
\widehat B=\twomat{0}{B}{B^*}{0},
\]
then, for all $k=1,\ldots,n$, we have
$\widehat\alpha(k) = \alpha(k)$ and $\widehat\alpha(2n+1-k) = -\alpha(k)$, and similar identities for
$\widehat\beta(k)$ and $\widehat\gamma(k)$.
Using these identities, (\ref{eq:TFW2}) reduces to the singular value inequality (\ref{eqh:5}), which ends the proof.
\qed
\section*{Acknowledgments}
We acknowledge support by an Odysseus grant from the Flemish FWO.
We are indebted to Weimin Miao for introducing us to the topic of this paper.


\begin{thebibliography}{9}
\bibitem{AK} K.M.R.~Audenaert and F.~Kittaneh, ``Problems and Conjectures in Matrix and Operator Inequalities'',
in: Operator Theory, J.~Zemanek ed., Banach Center Publications Series. (In Press). See also eprint arXiv:1201.5232.

\bibitem{bhatia} R.~Bhatia, \textit{Matrix Analysis}, Springer, Heidelberg (1997).

\bibitem{bhat} R.~Bhatia, ``Linear Algebra to Quantum Cohomology: The story of Alfred Horn's inequalities'',
The American Mathematical Monthly \textbf{108}(4), 289--318 (2001).

\bibitem{BU} J.-C.~Bourin and M.~Uchiyama, ``A matrix subadditivity inequality for $f(A+B)$ and $f(A)+f(B)$'',
Linear Algebra Appl.\ \textbf{423}(2--3), 512--518 (2007).



\bibitem{M} L.~Mirsky, ``Symmetric gauge functions and unitarily invariant norms'', 
Quart.\ J.\ Math.\ Oxford (2), \textbf{11}, 50--59 (1960).

\bibitem{thompson} R.C.~Thompson, ``Convex and concave functions of singular values of matrix sums'',
Pacific J.\ Math.\ \textbf{66}(1), 285--290 (1976).

\bibitem{TF} R.C.~Thompson and L.J.~Freede, ``On the eigenvalues of sums of Hermitian matrices'',
Linear Algebra Appl.\ \textbf{4}, 369--376 (1971).

\bibitem{uchiyama} M.~Uchiyama, Proc.\ Amer.\ Math.\ Soc.\ \textbf{134}(5), 1405--1412 (2005).

\bibitem{YS} M-C.~Yue and A.~M-C.~So, ``A Perturbation Inequality for the Schatten $p$-Quasi-Norm
and Its Applications in Low-Rank Matrix Recovery'', eprint arXiv 1209.0377v4 (2014).

\bibitem{ZQ} Y.~Zhang and L.~Qiu, ``From subadditive inequalities of singular values to triangle inequalities of canonical angles'',
SIAM J.\ Matrix Anal.\ Appl.\ \textbf{31}(4), 1606--1620 (2010).
\end{thebibliography}
\end{document}